\documentclass{mrlart2e}
\usepackage{amsmath, amsthm, url}
\usepackage{amssymb}

\newtheorem{thm}{Theorem}[section]
\newtheorem{prop}[thm]{Proposition} \newtheorem{cor}[thm]{Corollary}
\newtheorem{lemma}[thm]{Lemma} 
\newtheorem{lemmadef}[thm]{Lemma-Definition}

\theoremstyle{definition}
\newtheorem{defe}[thm]{Definition}

\theoremstyle{remark}
\newtheorem{exam}[thm]{Example}

\newtheorem{ntn}[thm]{Notation}
\newtheorem{ov}[thm]{Overview}
\newtheorem{rem}[thm]{Remark}
\numberwithin{equation}{section}

\def\GGS{\mathrm{GGS}}
\def\pr{\mathrm{pr}}
\def\peq{\underline{\prec}}

\def\O{\text{Ord}}

\def\Z{\mathbb Z}
\def\C{\mathbb C}
\def\Res{\text{Res}}
\def\o{\otimes}

\def\h{{\mathfrak h}}
\def\g{{\mathfrak g}}

\begin{document}
\title{Trigonometric solutions of the associative Yang-Baxter equation}
\author{Travis Schedler}

\begin{abstract}
We classify trigonometric solutions to the associative Yang-Baxter
equation (AYBE) for $A = \mathrm{Mat}_n$, the associative algebra of
$n$-by-$n$ matrices.  The AYBE was first presented in a 2000 article
by Marcelo Aguiar and also independently by Alexandre Polishchuk.
Trigonometric AYBE solutions limit to solutions of the classical
Yang-Baxter equation.  We find that such solutions of the AYBE are
equal to special solutions of the quantum Yang-Baxter equation (QYBE)
classified by Gerstenhaber, Giaquinto, and Schack (GGS), divided by a
factor of $q - q^{-1}$, where $q$ is the deformation parameter $q =
e^{\hbar}$.  In other words, when it exists, the associative lift of
the classical $r$-matrix coincides with the quantum lift up to a
factor. We give explicit conditions under which the associative lift
exists, in terms of the combinatorial classification of classical
$r$-matrices through Belavin-Drinfeld triples.  The results of this
paper illustrate nontrivial connections between the AYBE and both
classical (Lie) and quantum bialgebras.
\end{abstract}

\section{Introduction}
Let $A$ be an associative algebra (not necessarily with unit), and 
$r\in A\otimes A$. The {\it associative Yang-Baxter equation}
(AYBE) for $r$ over $A$ is the equation 
\begin{equation}
r^{13}r^{12}-r^{12}r^{23}+r^{23}r^{13}=0.
\end{equation}
This equation was introduced in \cite{Ag1,Ag2} and independently 
in \cite{P}. 

The algebraic meaning of this equation, explained 
in \cite{Ag1,Ag2}, is as follows. 
An associative algebra $A$ 
is called an infinitesimal bialgebra if it is equipped with 
a coassociative coproduct which is a derivation, i.e.~$\Delta(ab)
=(a\otimes 1)\Delta(b)+\Delta(a)(1\otimes b)$. 
This notion was introduced by Joni and Rota \cite{JR}
and is useful in combinatorics. 
Now, given an associative algebra $A$ and a 
solution $r\in A\otimes A$ of the AYBE, one can define 
a comultiplication by $\Delta(a)=(a\otimes 1)r-r(1\otimes a)$. 
(This comultiplication is a derivation for any $r$, and 
is coassociative if $r$ satisfies the AYBE). Thus, $(A,\Delta)$ is 
an infinitesimal bialgebra.  

One may also consider the AYBE with spectral parameter,
\begin{multline}
r^{13}(v_1-v_3)r^{12}(v_1-v_2)-r^{12}(v_1-v_2)r^{23}(v_2-v_3)+ \\
r^{23}(v_2-v_3)r^{13}(v_1-v_3)=0,
\end{multline}
where $r(v)$ is a meromorphic function of a complex variable $v$ with
values in $A \otimes A$.  Similarly to the usual (classical and
quantum) YBE, this is essentially the same equation, since $r(v)$ is a
solution of this equation if and only if $r(v \o 1 -1 \o v)$ satisfies the
usual AYBE over $\hat A$, where $\hat
A=A((v))$ is the algebra of Laurent series with coefficients in $A$, and
the tensor products $\hat A \o \hat A (\o \hat A)$ are completed in some form.

Further, one may consider a graded version of the AYBE. 
Namely, given a finite abelian group $\Gamma$, one may consider 
solutions $r$ of the usual AYBE over the
$\Gamma$-graded algebra $A \o \C[\Gamma]$
which have total degree zero, i.e.~are
sums of terms of bidegrees $(x,-x)\in \Gamma^2$.  
In this case, writing $r(u)$ for the part of $r$ of bidegree
$(u,-u)$, we obtain the following
equation for $r(u)$:
\begin{equation}
r^{13}(u+u')r^{12}(-u')-r^{12}(u)r^{23}(u+u')+r^{23}(u')r^{13}(u)=0.
\end{equation}

This equation, which one may call
the graded AYBE,
obviously makes sense for infinite groups $\Gamma$ as well;
moreover, if $\Gamma$ is a complex vector space, 
then one may require $r(u)$ to be meromorphic in $u$. 
Finally, as before, one can add a spectral parameter. 
In this form, 
(with a 1-dimensional space $\Gamma$),
the AYBE arose in the work of Polishchuk \cite{P}, in the study of 
$A_\infty$-categories attached to algebraic curves
of arithmetic genus 1. 
More precisely, the equation considered in 
\cite{P} is the graded AYBE with spectral parameter 
over the algebra 
$A^{op}$ opposite to $A$.  Using ordinary multiplication and
making the substitution $v = v_1-v_2$ and $v' = v_2-v_3$, the equation takes the form
\begin{multline} \label{aybe}
r^{12}(-u',v) r^{13}(u+u',v+v')
-r^{23}(u+u',v')r^{12}(u,v) \\ +r^{13}(u,v+v')r^{23}(u',v')=0,
\end{multline}
where $r$ is a meromorphic function of two complex variables with values in
$A\otimes A$. From now on, the term ``AYBE'' will be reserved for 
this equation. 

One special case studied in \cite{P} 
is where $A=\mathrm{Mat}_n(\C)$ and AYBE solutions $r(u,v)$ also satisfy
the unitarity condition
\begin{equation} \label{aun}
r^{21}(-u, -v) = -r(u,v),
\end{equation} and have a Laurent expansion near $u=0$ of the form
\begin{equation}
\label{laur}
r(u,v) = \frac{1 \o 1}{u} + r_0(v) + u r_1(v) + O(u^2).
\end{equation}
In this case, we will show that $r_0(v)$ satisfies the CYBE with spectral parameter,
\begin{equation} \label{cybepint}
[r_0(v)^{12}, r_0(v+v')^{13}] + [r_0(v)^{12}, r_0(v')^{23}] + [r_0(v+v')^{13}, r_0(v')^{23}] = 0,
\end{equation}
and the unitarity condition,
\begin{equation}
r(-v)^{21} = -r(v).
\end{equation}
This follows from the proof of the fact in \cite{P} that, even without
the Laurent condition \eqref{laur}, when the limit $\overline{r}(v) =
(\pr \o \pr)r(u,v) \bigl|_{u = 0}$ exists ($\pr$ is the projection away from
the identity to traceless matrices), it is a unitary solution
of the CYBE with spectral parameter.

In this paper we will classify all such matrices $r(u,v)$ where
$r_0(v) = \frac{\tilde r + e^v \tilde r^{21}}{1 - e^v}$ for $\tilde r$
a constant solution of the CYBE \eqref{cybepint} satisfying $\tilde r
+ \tilde r^{21} = \sum_{i,j} e_{ij} \o e_{ji}$.  These $\tilde r$ were
classified by Belavin and Drinfeld in the 1980's \cite{BD} in terms of
combinatorial objects known as Belavin-Drinfeld triples.  We will
discover that such matrices $r(u,v)$ correspond not to all
Belavin-Drinfeld triples for $\tilde r$, but to a subclass of them, called {\it
associative BD triples}.  In particular, we answer negatively the
question asked in Remark 1 of Section 5 of \cite{P}: whether any
nondegenerate solution $\overline{r}(v) = (\pr \o \pr) r_0(v)$ of the
CYBE can be ``lifted'' to such an AYBE solution $r(u,v)$ (see Remark
\ref{rbarrem}).  Also, for those triples which are associative, only
special classical $r$-matrices from the usual continuous family are
liftable.  Recall that the Belavin-Drinfeld classification assigns
to each BD triple a family of classical $r$-matrices parameterized by
a finite-dimensional vector space of skew-symmetric diagonal
components.  We will demonstrate that there is only a finite number
of choices of this component, up to scalars ($1 \o A +
A \o 1$), which yield an $r$-matrix liftable to an
associative $r$-matrix (this number is nonzero iff the BD triple is
associative).

More precisely, the condition for a classical $r$-matrix to be
``liftable'' to an associative $r$-matrix (a unitary solution of the
AYBE) satisfying \eqref{laur} is that the map $T: \Gamma_1 \rightarrow
\Gamma_2$ which defines the BD triple be ``liftable'' to a cyclic
permutation $\tilde T$ of the set $\{e_1, \ldots, e_n\}$. Here,
``liftable'' means that $T(\alpha_i) = \alpha_j$ implies that $\tilde
T(e_i) = \tilde T(e_j)$ and $\tilde T(e_{i+1}) = \tilde T(e_{j+1})$.
Here $\Gamma_1, \Gamma_2 \subset \Gamma= \{\alpha_1, \ldots,
\alpha_{n-1}\}$.
Then, the
skew-symmetric diagonal component $s$ parameterizing the $r$-matrices
for the BD triple is determined up to scalars by an explicit formula
(the solution to ``associative''
versions of the equations for $s$ in the CYBE theory).  There
are evidently finitely many choices of the lift $\tilde T$ of $T$, and
we define the BD triple to be ``associative'' if there exists
at least one.

We discover that such an associative $r$-matrix lifting a
classical $r$-matrix $\tilde r$ is closely
 related to the Gerstenhaber-Giaquinto-Schack 
(GGS) quantization of $\tilde r$, i.e.~a special  
matrix $R_\GGS(u) = 1 + u\tilde r +O(u^2)$ which satisfies the QYBE,
\begin{equation}
R^{12} R^{13} R^{23} = R^{23} R^{13} R^{12},
\end{equation}
and the Hecke condition,
\begin{equation}
(PR - q)(PR + q^{-1}) = 0,
\end{equation}
where $q = e^{u/2}$ and $P = \sum_{i,j} e_{ij} \o e_{ji}$ is the permutation
matrix (see \cite{GGS},\cite{S}). Namely, 
$R_{GGS}(u)=(e^{u/2} - e^{-u/2}) \lim_{v\to -\infty}r(u,v) = (e^{u/2} - e^{-u/2}) [r(u,v) - \frac{e^v}{1-e^v} P]$, where in the limit $v$ is taken to be real.

In fact, we can make the connection between the AYBE solution
``lifting'' a classical $r$-matrix (a solution of \eqref{cybe} and
\eqref{nu}) and the QYBE solution ``quantizing'' the classical
$r$-matrix more apparent by adding the spectral parameter $v$ back
into the quantum $R$-matrix.  That is, for any matrix $R = 1 + u r + O(u^2)$ satisfying
the QYBE and the Hecke condition which is a function only of the
parameter $q = e^{u/2}$, one can consider its ``Baxterization,''
\begin{equation}\label{rwsp}
R_{\mathrm{B}}(q,v) = \frac{e^v}{1 - e^v} (q - q^{-1}) P + R(q),
\end{equation}
which is a solution of the QYBE and the Hecke condition which quantizes the
CYBE solution with spectral parameter $\frac{r + e^v r^{21}}{1 - e^v}$ \cite{Mu}. 

Now, letting $R_{\mathrm{BGGS}}(q,v)$ be given by \eqref{rwsp} from $R_\GGS(q)$, we find that
\begin{equation} \label{rRc}
r(u,v) = \frac{R_{\mathrm{BGGS}}(q,v)}{q - q^{-1}}, \quad \mathrm{where\ } q = e^{u/2}.
\end{equation}
In particular,
this implies that the matrix $r(u,v)$ satisfies not only the AYBE but
also the QYBE.

\begin{rem}
The fact that the ``associative r-matrix'' $r(u,v)$ 
specializes to both classical and quantum r-matrices is 
in good agreement with the remark in \cite{Ag1} (p.2) 
that infinitesimal bialgebras
have nontrivial analogies and connections 
with both classical (Lie) and quantum bialgebras.
At the same time, we must admit that we don't have a conceptual explanation
for the validity of \eqref{rRc}.  To find such an explanation seems to be
an interesting problem.
\end{rem}

\begin{rem}
In \cite{Mu} (p.9) Mudrov quantizes certain Belavin-Drinfeld triples that
obey a slightly more restrictive version of the associative conditions than
those considered in this paper.  To do this, Mudrov uses the language
of associative Manin triples.  
It appears that the theory of \cite{Mu} is parallel
to \cite{Ag1,Ag2} and closely related to the content of this paper.
\end{rem} 

\begin{rem}
We expect that the results of this paper can be generalized to the
case of all trigonometric solutions of the CYBE with spectral parameter
(not just those obtained from constant CYBE solutions).  In this case,
we expect again that the classical $r$-matrices with spectral parameter
can be lifted provided they satisfy the BD associativity conditions and
the classical $r$-matrix for the triple is chosen correctly (in an
analogous way to the case of constant $r$-matrices).
Furthermore, for any given $r_0(v)$, the associative lift $r(u,v)$
should again be related to the quantum lift $R(q,v)$ by
\begin{equation}
R(q,v) = (q - q^{-1})r(u,v), \quad q=e^{u/2}.
\end{equation}
The matrix $R(q,v)$ should be given explicitly by a generalization of
the GGS formula (there is already a different kind of explicit formula
for $R(q,v)$ given in \cite{ESS} and \cite{ES}).
\end{rem}

\subsection{Acknowledgements}
I would like to thank Pavel Etingof for advising
and help with writing the introduction, 
and Alexander Polishchuk for posing the problem.  I am
grateful to the Harvard College Research Program for partially funding
this research.

\section{Background}
\begin{ov} We formally introduce Belavin-Drinfeld triples, the AYBE as presented in \cite{P}, and the GGS Conjecture \cite{GGS}, proved in \cite{S}.
\end{ov}

\subsection{Belavin-Drinfeld triples} \label{bd}

Let $(e_i), 1 \leq i \leq n,$ be the standard basis for $\mathbb C^n$.
Set $\Gamma = \{e_i - e_{i+1}: 1 \leq i \leq n-1\}$.
We will use the notation
$\alpha_i := e_i - e_{i+1}$.  Let $( , )$ denote the inner product
on $\mathbb C^n$ having $(e_i)$ as an orthonormal basis.

\begin{defe} \cite{BD}
A Belavin-Drinfeld triple of type $A_{n-1}$ is given by $(T, \Gamma_1,
\Gamma_2)$ where $\Gamma_1, \Gamma_2 \subset \Gamma$ and $T:
\Gamma_1 \rightarrow \Gamma_2$ is a bijection, satisfying

(a) $T$ preserves the inner product: $\forall \alpha, \beta
\in \Gamma_1$, $(T \alpha,T \beta) = (\alpha, \beta)$.

(b) $T$ is nilpotent: $\forall \alpha \in \Gamma_1, \exists k
\in \mathbb N$ such that $T^k \alpha \notin \Gamma_1$.
\end{defe}

Let $\mathfrak g = \mathfrak{gl}(n)$ be the Lie algebra of complex $n
\times n$ matrices. Define $\mathfrak h \subset \mathfrak g$ to be the
abelian subalgebra of diagonal matrices and $\mathfrak{g}' \subset \g$
to be the simple subalgebra of traceless matrices
(i.e.~$\mathfrak{sl}(n)$).  Elements of $\mathbb C^n$ define linear
functions on $\mathfrak h$ by $\bigl ( \sum_i \lambda_i e_i \bigr)
\bigl( \sum_i a_i\: e_{ii} \bigr)= \sum_i \lambda_i a_i$.  Under this
identification, we use $\Gamma$ as the set of simple roots of
$\mathfrak{g}'$ with respect to the Cartan subalgebra $\mathfrak{h} \cap \mathfrak{g}'$.
Let $P = \sum_{1 \leq i,j \leq n} e_{ij} \otimes e_{ji}$ be the
Casimir element inverse to the standard form, $(B,C) =
\mathrm{tr}(BC)$, on $\mathfrak g$.  It is easy to see that $P (w \o
v) = v \o w$, for any $v,w \in \C^n$, so we also call $P$ the permutation
matrix.  Let $P^0=\sum_i e_{ii} \o e_{ii}$ be the projection of $P$ to
$\mathfrak h \otimes \mathfrak h$.

For any Belavin-Drinfeld triple, consider the following equations for
$s \in \mathfrak h \wedge \mathfrak h$:

\begin{gather}
\label{tr02} \forall \alpha \in \Gamma_1,
\bigl[(\alpha - T \alpha) \otimes 1 \bigr] s = \frac{1}{2}
[(\alpha + T \alpha) \o 1] P^0.
\end{gather}

Belavin and Drinfeld showed that solutions $r \in \mathfrak{g} \o
\mathfrak g$ of the constant CYBE,
\begin{equation} \label{cybe}
[r^{12}, r^{13}] + [r^{12}, r^{23}] + [r^{13}, r^{23}] = 0,
\end{equation} satisfying
\begin{equation} \label{nu}
r + r^{21} = P = \sum_{i,j} e_{ij} \o e_{ji},
\end{equation} are given, up to inner isomorphism, by a discrete
datum (the Belavin-Drinfeld triple) and a continuous datum (a solution
$s \in \h \wedge \h$ of
\eqref{tr02}).  We now describe this classification.

For $\alpha = e_i - e_j$, set $e_\alpha := e_{ij}$.  Define
$|\alpha| = |j - i|$.  For any $Y \subset \Gamma$, set $\tilde Y =
\{\alpha \in \text{Span}(Y) \mid \alpha = e_i - e_j, i < j\}$ (the set
of positive roots of the semisimple subalgebra of $\mathfrak{g}'$ having $Y$
as its set of simple roots). In particular we will use the
notation $\tilde \Gamma, \tilde \Gamma_1$, and $\tilde \Gamma_2$.  We extend
$T$ additively to a map $\tilde \Gamma_1 \rightarrow \tilde \Gamma_2$,
i.e.~by $T(\alpha+\beta)=T \alpha +T \beta$. Whenever $T^k \alpha =
\beta$ for $k \geq 1$, we say $\alpha \prec \beta$ and
$O(\alpha,\beta) = k$, while $O(\beta, \alpha) = -k$.  Clearly $\prec$
is a partial ordering on $\tilde \Gamma$.  We will also use $\alpha
\peq \beta$ to denote that either $\alpha \prec \beta$ or $\alpha = \beta$.
Suppose $T^k \alpha = \beta$ for $\alpha = e_i - e_j$ and $\beta = e_l
- e_m$.  Then there are two possibilities on how $T^k$ sends $\alpha$
to $\beta$, since $T^k$ induces an isomorphism of the segment of the
Dynkin diagram corresponding to $\alpha$ onto the segment corresponding
to $\beta$.  Namely, either $T^k(\alpha_i) = \alpha_l$ and
$T^k(\alpha_{j-1})=\alpha_{m-1}$, or $T^k(\alpha_i) = \alpha_{m-1}$
and $T^k(\alpha_{j-1}) = \alpha_l$.  In the former case, call $T^k$
{\sl orientation-preserving on $\alpha$}, and in the latter, {\sl
orientation-reversing on $\alpha$}.  Let
\begin{equation}
C_{\alpha,\beta} =
\begin{cases} 1, & \text{if $T^k$ reverses orientation on $\alpha$,}
\\ 0, & \text{if $T^k$ preserves orientation on $\alpha$.}
\end{cases}
\end{equation}
Now, we define
\begin{gather} \label{adef}
a = \sum_{\alpha \prec \beta} (-1)^{C_{\alpha,\beta} (|\alpha|-1)}
(e_{-\alpha} \o e_\beta - e_{\beta} \o e_{-\alpha}),
\\ \label{rdef} r_{st} = \frac{1}{2} \sum_i e_{ii} \o e_{ii} +
\sum_{\alpha \in \tilde \Gamma} e_{-\alpha} \o e_{\alpha},
\quad r_{T, s} = s + a + r_{st}.
\end{gather}
Here $r_{st} \in \mathfrak{g} \o \mathfrak{g}$ is the standard solution of
the CYBE satisfying $r_{st} + r_{st}^{21} = P$, and $r_{T, s}$ is
the solution of the CYBE corresponding to the data
$((\Gamma_1,\Gamma_2,T), s)$ ($r_{st}$ corresponds to the trivial BD triple with $s = 0$).  It follows from \cite{BD} that
\begin{prop}\cite{BD}
Any solution $\tilde r \in \mathfrak{g}$ of \eqref{cybe} and
\eqref{nu} is equivalent to a solution $r_{T,s}$ given in \eqref{rdef} for
some Belavin-Drinfeld triple and continuous parameter $s$, under an
inner automorphism of $\mathfrak{g}$.
\end{prop}
\begin{defe} Solutions of \eqref{cybe} and \eqref{nu} will be called
{\sl classical $r$-matrices}.
\end{defe}

\begin{exam} \label{gcge} For a given $n$, there
are exactly $\phi(n)$ BD triples ($\phi$ is the Euler $\phi$-function)
in which $|\Gamma_1| + 1 = |\Gamma|$ \cite{GG}.  These are called {\sl
generalized Cremmer-Gervais} triples (the usual Cremmer-Gervais triple
is the 
special case $m=1$ in the following classification).  These are
indexed by $\{m \in \{1, \ldots, n\} \mid \text{gcd}(n,m) = 1 \}$, and
given by $\Gamma_1 = \Gamma \setminus \{\alpha_{n-m}\}$, $\Gamma_2 =
\Gamma \setminus \{\alpha_m\}$, and $T(\alpha_i) =
\alpha_{\text{Res}(i+m)}$, where $\Res$ gives the residue modulo $n$
in $\{1,\ldots,n\}$.  For these triples, there is a unique $s$ taken
to lie in $\mathfrak{g}' \wedge \mathfrak{g}'$, given by
$s^{ii}_{ii} = 0, \forall i$, and $s_{ij}^{ij} = \frac{1}{2} -
\frac{1}{n}\text{Res}(\frac{j-i}{m}), i \neq j$ (this is easy to
verify directly and is also given in \cite{GG}). We will see that this
formula for $s$ generalizes to formula \eqref{stp} in the associative
case.
\end{exam}

\subsection{The CYBE and AYBE with parameters}
The CYBE takes the following form ``with spectral parameter'' over a Lie
algebra $\mathfrak{a}$:
\begin{equation} \label{cybep}
[r^{12}(x), r^{13}(x+y)] + [r^{12}(x), r^{23}(y)] + [r^{13}(x+y),
 r^{23}(y)] = 0.
\end{equation}
Here $r(v)$ is a meromorphic function of $v$ with values in $\mathfrak{a} \o
\mathfrak{a}$.  A solution $r$ is called {\sl unitary} if
\begin{equation} \label{cun}
r(v) = -r^{21}(-v).
\end{equation}
\begin{lemmadef}  If $r$ is a constant solution of the CYBE, then
$\frac{r + e^v r^{21}}{1 - e^v}$ is a unitary solution of the CYBE
with spectral parameter $v$. For any constant solution $r$, define
\begin{equation} \label{adsp}
\hat r(v) = \frac{r + e^v r^{21}}{1 - e^v}.
\end{equation}
\end{lemmadef}
\begin{proof} This follows immediately. \end{proof}

The version of the AYBE we consider has the form
\begin{multline}
r^{12}(-u',v) r^{13}(u+u',v+v') - r^{23}(u+u',v') r^{12}(u,v) \\ +
r^{13}(u,v+v') r^{23}(u',v') = 0.
\end{multline}
The {\sl unitarity} condition is
\begin{equation}
r^{21}(-u,-v) = -r(u,v).
\end{equation}
Unitary AYBE solutions give rise to CYBE solutions in the following way:
\begin{prop} \cite{P}\label{ac}
Let $A = \mathfrak{g}$ and $\pr: \mathfrak{g} \rightarrow
\mathfrak{g}'$ the orthogonal projection with respect to the standard
form,
$(B,C) = \mathrm{tr}(BC)$.  If $r(u,v)$ is a unitary solution of the
AYBE, and the limit $\overline{r}(v) = [(\pr \o \pr) r(u,v)]|_{u=0}$
exists, then $\overline{r}(v)$ is a unitary solution of the CYBE with
spectral parameter.
\end{prop}
\begin{proof} 
We repeat the proof of \cite{P} (since it is short and
we will use \eqref{acy} later). 
First note that the unitarity of $\overline{r}$ follows immediately
from the unitarity of $r$.
Substituting $r^{21}(-u,-v) = -r^{12}(u,v)$, we rewrite the AYBE as
\begin{multline}
-r^{21}(u', -v) r^{13}(u+u', v+v') + r^{23}(u+u', v') r^{21}(-u,-v) \\ +
r^{13}(u,v+v') r^{23}(u', v') = 0.
\end{multline}
We permute the first two components, yielding
\begin{multline}
-r^{12}(u', -v) r^{23}(u+u', v+v') + r^{13}(u+u', v') r^{12}(-u, -v)
\\ + r^{23}(u, v+v') r^{13}(u', v') = 0.
\end{multline}
This resembles the AYBE with the order of each product reversed (which
we seek). To obtain it, we make the linear change of 
variables given by $u \mapsto u', u' \mapsto u, v \mapsto -v$, and 
$v' \mapsto v + v'$:
\begin{multline}
r^{13}(u+u', v+v') r^{12}(-u', v) - r^{12}(u, v)r^{23}(u+u', v') 
\\ + r^{23}(u', v') r^{13}(u, v+v') = 0.
\end{multline}
Subtracting this from the AYBE, we get
\begin{multline} \label{acy}
[r^{12}(-u', v), r^{13}(u+u', v+v')] + [r^{12}(u,v), r^{23}(u+u', v')]
\\ + [r^{13}(u, v+v'), r^{23}(u', v')] = 0.
\end{multline}
Applying $\pr \o \pr \o \pr$, we get the same equation with $(\pr \o \pr) r$
replacing $r$, and then we may take the limit $u \rightarrow 0$ to find
that $\overline{r}(v)$ satisfies the CYBE with spectral parameter.
\end{proof}
This warrants the
\begin{defe}
Solutions of \eqref{aybe} and \eqref{aun} are called {\sl associative
$r$-matrices}.
\end{defe}

In the case we consider, $r(u,v)$ has a Laurent expansion at $u=0$ of the
form \eqref{laur}, and this result can be strengthened:
\begin{lemma}
If $r(u,v)$ has a Laurent expansion at $u=0$ of the form 
\begin{equation*}
r(u,v) = \frac{1 \o 1}{u} + r_0(v) + u r_1(v) + O(u^2)
\end{equation*} 
and is an associative $r$-matrix, then $r_0(v)$ is a solution of the CYBE
with spectral parameter.
\end{lemma}
\begin{proof}
This follows from \eqref{acy}, since $1$ commutes with anything.
\end{proof}

\subsection{The GGS quantization}
Given any Belavin-Drinfeld triple $(\Gamma, \Gamma', T)$ and
any matrix $s \in \h \wedge \h$ satisfying \eqref{tr02}, the CYBE solution
$r_{T, s}$ is one-half the linear term in $\hbar$ of a quantum $R$-matrix $R_\GGS = 1 + 2 r_{T, s} \hbar + O(\hbar^2)$, which satisfies the quantum Yang-Baxter equation,
\begin{equation} \label{qybe}
R^{12} R^{13} R^{23} = R^{23} R^{13} R^{12},
\end{equation}
and the Hecke relation,
\begin{equation}\label{hecke}
(PR - q)(PR + q^{-1}) = 0, \quad q = e^\hbar.
\end{equation}
The matrix $R_\GGS$ is given by a simple (yet not fully understood)
formula proposed by Gerstenhaber, Giaquinto, and Schack \cite{GGS} in
1993:
\begin{gather}
R_{st} = q \sum_{i} e_{ii} \o e_{ii} + \sum_{i \neq j} e_{ii} \o
e_{jj} + (q - q^{-1}) \sum_{\alpha > 0} e_{-\alpha} \o e_{\alpha}, \\
R_\GGS = q^s \biggl(R_{st} + (q - q^{-1}) \sum_{\alpha \prec \beta}
(-1)^{C_{\alpha, \beta} (|\alpha| - 1)} [q^{-C_{\alpha,
\beta}(|\alpha| - 1) - \mathrm{PS}(\alpha, \beta)} e_{-\alpha} \o
e_{\beta} \nonumber \\ - q^{C_{\alpha, \beta}(|\alpha| - 1) +
\mathrm{PS}(\alpha, \beta)} e_{\beta} \o e_{-\alpha}]\biggr) q^s, \label{ggsf}
\end{gather}
where $\mathrm{PS}(\alpha, \beta)$ is defined as follows.  First, we
define the relation $\alpha \lessdot \beta$ for $\alpha > 0, \beta >
0$ to mean that, writing $\alpha = e_i - e_j$ and $\beta = e_k - e_l$,
we have $j = k$.  In other words, considering $\alpha$ to be the line
segment with endpoints $i$ and $j$ and $\beta$ the line segment with
endpoints $k$ and $l$ on the real line, we have that $\alpha$ lies
adjacent to $\beta$ on the left.  Now, let $[\mathrm{statement}] = 1$
if ``statement'' is true, and $[\mathrm{statement}] = 0$ otherwise.
Then, $\mathrm{PS}$ is given by
\begin{multline}
\mathrm{PS}(\alpha, \beta) = \frac{1}{2} \bigl([\alpha \lessdot \beta]
+ [\beta \lessdot \alpha]\bigr) + [\exists \gamma \mid \alpha \prec
\gamma \prec \beta, \alpha \lessdot \gamma] \\ + [\exists \gamma \mid
\alpha \prec \gamma \prec \beta, \gamma \lessdot \alpha].
\end{multline}
\begin{thm}[The GGS Conjecture] \cite{GGS}, \cite{S} The element $R_\GGS$ satisfies
the QYBE \eqref{qybe} and the Hecke condition \eqref{hecke}.
\end{thm}

\section{Statement of the main theorem}
\begin{ov}
In this section, we state the main theorem, which gives (1) the
associativity conditions under which a classical $r$-matrix can be
lifted to an associative $r$-matrix, (2) the formula relating the
associative $r$-matrix to the GGS quantum $R$-matrix, and (3) a new,
explicit formula for the GGS $R$-matrix in this case (which is a
generalization of Giaquinto's formula for the GGS $R$-matrix in the
case of generalized Belavin-Drinfeld triples).
\end{ov}

\begin{defe} \label{atd}
Call a triple an {\sl associative triple} if (i) the triple preserves
orientation, and (ii) there exists a cyclic permutation $\tilde T$ of
$\{1,\ldots,n\}$ such that $T(\alpha_i) = \alpha_j$ implies $\tilde
T(i) = j$ and $\tilde T(i+1) = j+1$. Such permutations are called {\sl
compatible permutations}.  The structure $(\Gamma_1, \Gamma_2, T,
\tilde T)$ is called an {\sl associative structure}.  Given such a
structure, we define for each $i,j \in \{1, \ldots, n\}$ the function
$O(i,j)$ to be the least nonnegative integer such that $\tilde
T^{O(i,j)}(i) = j$.
\end{defe}

\begin{ntn} We will use the notation $s_0 = (\pr \o \pr) s$ in the future.
\end{ntn}

\begin{thm}\label{mt}
(1a) A classical $r$-matrix $\hat r_{T,s}$ 
is the zero-degree term $r_0(v)$ of the Laurent expansion
\eqref{laur} of an associative $r$-matrix $r(u,v)$ iff $(\Gamma_1,
\Gamma_2, T)$ is associative, and $s_0 = (\pr \o \pr) s$ is given by
the formula
\begin{equation} \label{stp}
s_0 = \sum_{i \neq j} \bigl( \frac{1}{2} - \frac{O(i,j)}{n} \bigr) e_{ii} \o
e_{jj},
\end{equation}

(1b) or equivalently satisfies
\begin{equation} \label{tra}
[(e_i - e_{\tilde T(i)}) \o 1] s_0 = \frac{1}{2} [(e_i + e_{\tilde T(i)}) \o 1]
[(\pr \o \pr) P^0].
\end{equation}

(2a) In this case, there is a unique associative $r$-matrix having a Laurent
expansion of the form
\begin{equation} \label{tle}
\frac{1 \o 1}{u} + \hat r_{T,s}(v) + O(u),
\end{equation}
and it is given by
\begin{equation} \label{afgq}
r(u,v) = \frac{e^v}{1 - e^v} P + \frac{R_{\mathrm{GGS}}(e^{u/2})}{e^{u/2} - e^{-u/2}},
\end{equation}
where $R_{\mathrm{GGS}}(e^{u/2})$ is the GGS matrix for the same $T$ and $s$ as $r_{T,s}$, replacing $q$ by $e^{u/2}$.

(2b) Using the Baxterization $R_{\text{BGGS}}(q,v)$, we get
\begin{equation} \label{bafgq}
r(u,v) = \frac{R_{\mathrm{BGGS}}(e^{u/2}, v)}{e^{u/2} - e^{-u/2}}.
\end{equation}

(3) The matrix $R_\GGS(q)$ occurring in \eqref{afgq} is given by
\begin{multline} \label{gqf}
R_{\mathrm{GGS}}(q) = q^{s - s_0} \biggl[ \sum_{i,j}
q^{1 - 2 O(i,j)/n} e_{ii} \o e_{jj} \\ + (q - q^{-1}) \biggl(
\sum_{\alpha > 0} e_{-\alpha} \o e_\alpha + \sum_{\alpha \prec \beta}
\bigl( q^{-2O(\alpha,\beta)/n} e_{-\alpha} \o e_\beta - q^{2O(\alpha,
\beta)/n} e_{\beta} \o e_{-\alpha} \bigr) \biggr) \biggr] q^{s-s_0}
\end{multline}
for any associative structure $(\Gamma_1, \Gamma_2,
T, \tilde T)$, where $s_0 = (\pr \o \pr) s$ is determined by \eqref{stp}.
\end{thm}

\begin{rem} \label{rbarrem}
One can also classify associative $r$-matrices where we require only
that the limit $\overline{r}(v) = (\pr \o \pr)(r(u,v)) \bigl|_{u=0}$
exist and satisfy $\overline{r}(v) = \hat {\tilde r}$
for some classical $r$-matrix $\tilde r$ over
$\mathfrak{g}'$. When the Laurent condition \eqref{laur} holds, 
all such lifts of $\overline{r}$ (without
fixing $r_0$) are
equal to $e^{cuv} r'(u,v)$, for $r'(u,v)$ an associative $r$-matrix
classified in Theorem \ref{mt} and $c \in \C$.  To see this, first
note that the BD associativity and $s_0$ conditions must still be
satisfied, because our proof of this part only uses the projection of
the AYBE away from scalars.  (This observation answers negatively the
question asked in Remark 1 of Section 5 in \cite{P}: whether, for any
unitary nondegenerate $\mathfrak{g}' \otimes \mathfrak{g}'$-valued
CYBE solution $\overline{r}(v)$ with spectral parameter, there exists
a unitary AYBE solution $r(u,v)$ having a Laurent expansion at $u=0$
of the form $r(u,v) = \frac{1\o 1}{u} + r_0(v) + O(u)$, such that
$(\pr \o \pr) r_0(v) = \overline{r}(v).$) Then, the result follows
from the fact (using Remark 2 in Section 5 of \cite{P}) that any two
associative $r$-matrices $r(u,v), r'(u,v)$ with Laurent expansions of
the form \eqref{laur} such that $\overline{r}(v) = \overline{r'}(v)$
are related by $r_0(v) - r'_0(v) = (1 \o 1)cv + \Phi^1 - \Phi^2$ where $c \in
\C$ and $\Phi \in \h$ satisfies $(\alpha, \Phi) = (T \alpha, \Phi),
\forall \alpha \in \Gamma_1$.  In this paper, we focus on lifts of
$r_0$ when it is a classical $r$-matrix, rather than lifting just
$\overline{r}$, since the result is cleaner.
\end{rem}

\begin{rem}
Equation \eqref{tra} can be thought of as the ``associative'' version
of \eqref{tr02} classifying classical $s$; it just so happens in the
associative case that these equations completely determine $s_0$ by
the choice of $\tilde T$.
\end{rem}

\begin{rem}
In the case of generalized Cremmer-Gervais triples (see Remark
\ref{gcge}), \eqref{gqf} is the formula found by Giaquinto \cite{S}.
Indeed, a generalized Cremmer-Gervais triple has a unique associative
structure, under which \eqref{stp} becomes the formula given in Remark
\ref{gcge}.
\end{rem}

\begin{rem}
Note that, given any associative choice of $T$, there are finitely
many possible compatible choices of $\tilde T$ (depending on $T$, and
up to $(n-1)!$ for the case of $T$ trivial).  Hence, the space of
associative matrices for each associative triple is parameterized by a
finite parameter ($\tilde T$) and a continuous parameter (the choice
of $s - s_0$).  The matrix $s-s_0$ can be any element in $\Lambda^2 \h
\cap (1 \o \h + \h \o 1)$ satisfying $[(\alpha - T \alpha) \o 1]
(s-s_0) = 0, \forall \alpha \in \Gamma_1$. In other words, $s - s_0 =
\Phi \o 1 - 1 \o \Phi$ for $\Phi \in \h$ any element satisfying
$(\alpha, \Phi) = (T \alpha, \Phi), \forall \alpha \in \Gamma_1$.
\end{rem}
\section{Proof of the main theorem \eqref{mt}}
\begin{ov} 
We prove the parts of Theorem \ref{mt} in the reverse order. Thus, in
the first subsection, we prove part (3), namely the explicit formula
for $R_\GGS$ for associative BD triples where $s_0 = (\pr \o \pr) s$
is given by \eqref{stp} for a choice of a compatible permutation
$\tilde T$.  Then, in the second subsection, we prove parts (2a) and
(2b) of Theorem \ref{mt}, namely verifying that $r(u,v) =
\frac{R_\GGS(e^{u/2})}{e^{u/2} - e^{-u/2}} + \frac{e^v}{1 - e^v} P$ in
fact satisfies the AYBE and unitarity conditions and lifts the
classical $r$-matrix, and is the unique such element.  Finally, in the
third subsection, we prove part (1) of Theorem \ref{mt}, that the BD
associativity and $s_0$-compatibility conditions are necessary and
sufficient for the lift to exist (necessity is all that will remain).
\end{ov}

\subsection{Proof of Theorem \ref{mt}, part (3): the generalization
of Giaquinto's formula}
\begin{ov}
We prove the generalization of Giaquinto's
formula \eqref{gqf} via a straightforward computation.
\end{ov}

First, we prove a lemma which gives a new formula for the combinatorial
constant $\mathrm{PS}(\alpha, \beta)$:
\begin{lemma}
For any $\alpha \prec \beta$, the number $\mathrm{PS}(\alpha, \beta) =
1 - (\alpha \o \beta) s$.
\end{lemma}
\begin{proof}
Note that, for $\beta = T^k \alpha$ ($k \geq 1$),
\begin{multline}
(\alpha \otimes \beta) s = \sum_{i = 0}^{k-1} [(T^i \alpha - T^{i+1} \alpha) \o \beta] s
= \sum_{i = 0}^{k-1} \frac{1}{2} (T^i \alpha + T^{i+1} \alpha, \beta) \\ = \frac{1}{2} (\alpha,
\beta) + \frac{1}{2} (2) + \sum_{i = 1}^{k-1} (T^i \alpha, \beta) \\
= 1 - \frac{1}{2} \bigl([\alpha \lessdot \beta] + [\beta \lessdot \alpha] \bigr) - \bigl([\exists \gamma
\mid \alpha \prec \gamma \prec \beta, \gamma \lessdot \beta] + [\exists \gamma \mid \alpha \prec \gamma \prec \beta, \beta \lessdot \gamma]\bigr),
\end{multline}
which proves the desired result.
\end{proof}

\begin{cor} The matrix $R_\GGS$ can be written as
\begin{multline} \label{nggs}
(q - q^{-1}) \biggl[\sum_{\alpha} e_{-\alpha} \o
e_\alpha  + \sum_{\alpha = e_i - e_j \prec \beta = e_k - e_l} (-1)^{C_{\alpha, \beta}(|\alpha| - 1)}
\\ \bigl(q^{-C_{\alpha, \beta}(|\alpha| - 1) + s_{ik}^{ik} + s_{jl}^{jl} - 1} e_{-\alpha} \o e_{\beta}  -
q^{C_{\alpha, \beta}(|\alpha| - 1) + 1 - s_{ik}^{ik} -
s_{jl}^{jl}} e_{\beta} \o e_{-\alpha}\bigr)\biggr]
+ q^{\sum_i e_{ii} \o e_{ii} + 2s}.
\end{multline}
\end{cor}
\begin{proof}
This follows immediately by expanding $(\alpha \o \beta) s = s_{ik}^{ik} + s_{jl}^{jl} - s_{il}^{il} - s_{jk}^{jk}$ for $\alpha = e_i - e_j$ and $\beta = e_k - e_l$, and noticing that $q^s e_{-\alpha} \o e_\beta q^s = q^{s_{jk}^{jk} + s_{il}^{il}}$ in this case.
\end{proof}
\begin{proof}[Proof of Theorem \ref{mt}, part (3).]
In the associative case where $s_0 = (\pr \o \pr)s$ is given by \eqref{stp} 
for a
compatible permutation $\tilde T$, we can simplify \eqref{nggs}. Let us
assume first that $s_0 = s \in \Lambda^2 \mathfrak{g}'$. 
Then, for each $\alpha = e_i - e_j \prec \beta = e_k - e_l$, we have
$C_{\alpha, \beta} = 0$ and $s_{ik}^{ik} = s_{jl}^{jl} = \frac{1}{2} -
\frac{O(i,k)}{n}$.  So, we rewrite
\eqref{nggs} as follows:
\begin{multline}
R_\GGS = \sum_{i,j} q^{1-2O(i,j)/n} e_{ii} \o e_{jj} \\ + (q-q^{-1})
\biggl[\sum_{\alpha > 0} e_{-\alpha} \o e_\alpha + \sum_{\alpha \prec
\beta} \bigl(q^{-2O(\alpha, \beta)/n} e_{-\alpha} \o e_\beta -
q^{2O(\alpha, \beta)/n} e_\beta \o e_{-\alpha}\bigr) \biggr].
\end{multline}

In the general case where $s$ is not necessarily equal to $s_0$, the
result follows from the fact, evident in \eqref{ggsf}, that $R_\GGS =
q^{s-s'} R' q^{s-s'}$, where $R'$ is the GGS matrix for the same
triple as $R_\GGS$, but replacing $s$ with $s'$.
\end{proof}

\subsection{Proof of Theorem \ref{mt}, parts (2a) and (2b): the GGS
$R$-matrix satisfies the AYBE with slight modifications}
\begin{ov} 
We verify that the $r(u,v)$ given by
\eqref{afgq} and \eqref{gqf} satisfies the AYBE and the unitarity
condition by a direct computation using BD combinatorics.
 A lemma from \cite{P} proves that $r(u,v)$
is uniquely determined by $r_0$ in \eqref{laur}, it is easy to check
that $r(u,v)$ lifts $r_{T,s}$ (i.e.~that $r_0 = r_{T,s}$).
These results prove part (2a) of
Theorem \ref{mt}, from which (2b) immediately follows.
As in the previous subsection, most of the
work reduces to the case where $s = s_0 \in \Lambda^2
\mathfrak{g}'$.
\end{ov}

\begin{lemma} \label{ruvc}
Fix some associative structure $(\Gamma_1, \Gamma_2, T, \tilde T)$
and choice of $s$ such that $s_0$ is given by \eqref{stp}.
Let $r(u,v)$ be given by \eqref{afgq}.  Let $r_0(v)$ be the
classical $r$-matrix which
is the term of degree-zero in the Laurent expansion of $r(u,v)$
in $u$ at $u=0$.  Then $r_0(v) = r_{T,s}$.
\end{lemma}

\begin{proof}
This follows from a simple computation using the next lemma \eqref{gruvl}.
Alternatively, it follows from the connection between $R_\GGS$ and
$r_{T,s}$.
\end{proof}

\begin{lemma} \label{gruvl} Set $s - s_0 = \Phi^1 - \Phi^2$ where
$\Phi \in \h$ satisfies $(\alpha, \Phi) = (T \alpha, \Phi)$ for any $\alpha \in \Gamma_1$.  Using
\eqref{gqf}, we can write the matrix $r(u,v)$ given by \eqref{afgq} as
follows:
\begin{multline}
r(u,v) = \frac{e^v}{1 - e^v} P + e^{-\Phi^2 u} \biggl[ \frac{1}{1
- e^{-u}} \sum_{i,j} e^{-O(i,j)u/n} e_{ii} \o e_{jj} \\ + \sum_{\alpha > 0}
e_{-\alpha} \o e_\alpha + \sum_{\alpha \prec \beta} \bigl(
e^{-\O(\alpha,\beta) u/n} e_{-\alpha} \o e_\beta - e^{\O(\alpha,\beta)
u/n} e_{\beta} \o e_{-\alpha} \bigr) \biggr] e^{\Phi^1 u}. \label{gruv}
\end{multline}
\end{lemma}
\begin{proof}
By the definition of $s_0$, we may write $s - s_0 = \Phi^1 - \Phi^2$.
The fact that $\Phi \in \h$ satisfies $(\alpha, \Phi) = (T \alpha,
\Phi)$ for any $\alpha \in \Gamma_1$ follows directly from the fact
that $[(\alpha - T \alpha) \o 1] (s - s_0) = 0$.  Now, it follows that
$e^{(\Phi^1 + \Phi^2) u}$, or simply $\Phi^1 + \Phi^2$, commutes with
$R_\GGS$.  Together with the fact that $e^t P e^t = 0$ for $t$ any
skew-symmetric matrix, we find that
\begin{equation}
r(u,v) = \frac{e^v}{1 - e^v} P + e^{-\Phi^2 u} \frac{R^0_\GGS(e^{u/2})}{e^{u/2} - e^{-u/2}} e^{\Phi^1 u},
\end{equation}
where $R^0_\GGS$ is the GGS matrix quantizing $r_{T,s_0}$.  Now \eqref{gruv}
follows from \eqref{gqf} with a small amount of manipulation. 
\end{proof}

\begin{ntn} For any $A \o A$-valued function $t$ of $u$ and $v$ (possibly constant
in one or both variables), we will denote by $AYBE(t)$ the LHS
of \eqref{aybe}.
\end{ntn}

\begin{lemma} \label{p12l}
Suppose that $r(u,v)$ is a solution of the AYBE and
$\Phi \in \mathfrak{h}$ is any diagonal matrix such that $\Phi^1 + \Phi^2$
commutes with $r(u,v)$.  The element
\begin{equation}\label{ors}
r'(u,v) = e^{-\Phi^2 u} r(u,v) e^{\Phi^1 u}
\end{equation}
also satisfies the AYBE. If, in addition, $r(u,v)$ is unitary, 
then so is $r'(u,v)$.
\end{lemma}
\begin{proof}
It is clear that $r'(u,v)$ satisfies the unitarity condition iff
$r(u,v)$ does. So, we show that $r'(u,v)$ satisfies the AYBE if
$r(u,v)$ does.  
Since $[\Phi^1 + \Phi^2, r] = 0$, it follows that $e^{(\Phi^1 + \Phi^2) z}$
commutes with $r$ for any complex variable $z$.  We make use of this
fact in the following computation, setting $t(u,v) = e^{-\Phi^2 u}
r(u,v) e^{\Phi^1 u}$:
\begin{multline}
t^{12}(-u', v)t^{13}(u+u', v+v') \\ = e^{\Phi^2 u'} r^{12}(-u', v) e^{-\Phi^1 u' - \Phi^3(u+u')}
r^{13}(u+u', v+v') e^{\Phi^1 (u+u')} \\ =
e^{\Phi^2 u' - \Phi^3 u} r^{12}(-u', v) r^{13}(u+u', v+v') e^{\Phi^1 u - \Phi^3 u'},
\end{multline}
and similarly,
\begin{gather}
t^{23}(u+u', v')t^{12}(u, v) = e^{\Phi^2 u' - \Phi^3 u} r^{23}(u+u',
v') r^{12}(u, v) e^{\Phi^1 u - \Phi^3 u'}, \mathrm{and}\\ t^{13}(u, v+v')
t^{23}(u', v') = e^{\Phi^2 u' - \Phi^3 u} r^{13}(u, v+v') r^{23}(u',
v') e^{\Phi^1 u - \Phi^3 u'}.
\end{gather}
Hence, it follows that $AYBE(t) = e^{\Phi^2 u' - \Phi^3 u} AYBE(r)
e^{\Phi^1 u - \Phi^3 u'}$. Hence, $t$ satisfies the AYBE iff $r$ does,
proving the desired result.
\end{proof}

\begin{lemma} \label{nypl}
The element $r(u,v) = y(u) + \frac{e^v}{1 - e^v} P$ satisfies the AYBE
and the unitarity condition, where $y(u)$ is any solution of the AYBE such that
$y(-u) + y^{21}(u) = P$. 
\end{lemma}
\begin{proof}  
Using facts of the form $P^{12} t^{13} = t^{23} P^{12}$ which follow
because $P$ is the permutation matrix, we compute
\begin{multline}
AYBE(y + f(v) P) = AYBE(y) + f(v+v') [y^{12}(-u') + y^{21}(u')] P^{13}
\\ + [f(v)f(v+v') - f(v')f(v) + f(v+v')f(v')]P^{12} P^{13}
\\ = [f(v+v')  +  f(v)f(v+v') - f(v') f(v) + f(v+v')f(v')] P^{12} P^{13}.
\end{multline}
So, the AYBE is satisfied for $y + f(v) P$, where $f$ is any function
satisfying the relation
\begin{equation}
f(v+v') = \frac{f(v)f(v')}{1 + f(v) + f(v')}.
\end{equation}
We can rewrite this as
\begin{equation}
f(v+v')^{-1} = f(v)^{-1} + f(v')^{-1} + f(v)^{-1} f(v')^{-1},
\end{equation}
which is the same as the condition that $g(v) = f(v)^{-1} + 1$
satisfies $g(v+v') = g(v) g(v')$.  So the solutions are $g(v) =
e^{Kv}$ for $K \in \C$, and in particular, when $K = -1$, we find
$f(v) = \frac{e^v}{1 - e^v}$.

Furthermore, provided $K \neq 0$, we evidently have $\frac{1}{e^{Kv} - 1}
+ \frac{1}{e^{-Kv} - 1} = -1$, so that $y(u) + \frac{e^{-Kv}}{1 - e^{-Kv}} P$
satisfies the unitarity condition.  This concludes the proof.
\end{proof}

\begin{rem}
This lemma essentially shows how to ``Baxterize'' AYBE solutions. As
mentioned in the introduction, we know that the same procedure
works for QYBE solutions using a result from \cite{Mu}.
\end{rem}

\begin{lemma}\label{ypl}
The element
\begin{multline}
y(u) = \frac{1}{1
- e^{-u}} \sum_{i,j} e^{-O(i,j)u/n} e_{ii} \o e_{jj} \\ + \sum_{\alpha > 0}
e_{-\alpha} \o e_\alpha + \sum_{\alpha \prec \beta} \bigl(
e^{-\O(\alpha,\beta) u/n} e_{-\alpha} \o e_\beta - e^{\O(\alpha,\beta)
u/n} e_{\beta} \o e_{-\alpha} \bigr) \label{yeq}
\end{multline}
satisfies the AYBE.
\end{lemma}
\begin{proof}
We will compute the coefficients $AYBE(y)_{ikm}^{jlp}$ and see that
they are all zero, so that $y$ satisfies the AYBE.  Note that we need
only check those indices for which $i+k+m = j+l+p$, because all
nonzero coefficients in the formula for $AYBE(y)$ obey this relation,
and the product or sum of matrices whose nonzero coefficients obey
this relation yields another matrix of the same form.

First, let us compute the coefficient $AYBE(y)_{ikm}^{jlp}$ for $i \neq
j, k \neq l,$ and $m \neq p$, subject to the relation $i+k+m = j+l+p$.
We have
\begin{multline} \label{aby}
AYBE(y)_{ikm}^{jlp} = y_{ik}^{\_l}(-u') y_{\_m}^{jp}(u+u') \\ -
y_{km}^{\_p}(u) y_{i\_}^{jl}(u+u') + y_{im}^{j\_}(u) y_{k\_}^{lp}(u'),
\end{multline}
where the underscore means that the index is deduced from the other
three by setting equal the sums of the upper and lower indices.  In
each product of two coefficients, the two underscores are equal.

We claim that either two or none of the three terms on the right-hand
side are nonzero, and that when there are two nonzero terms, they
cancel.  To see this, set $\alpha = e_i - e_j, \beta
= e_k - e_l$, and $\gamma = e_m - e_p$.  Suppose that $|\alpha| = |i-j| >
|\beta| = |k - l|$ and $|\alpha| > |\gamma|$.  Then if the first term
in the RHS of \eqref{aby} is nonzero, it follows that $-\alpha = T^{c}
\beta + T^{d} \gamma$, for some $c, d \in \Z$, and furthermore that
exactly one of the other two terms is nonzero: the second term if $|c|
< |d|$, the third term if $|d| < |c|$, or if $c = d$ then the second
term is nonzero iff $\alpha < 0$ (and the third term iff $\alpha >
0$).  Conversely, if the second or third term is nonzero, then the
first term must be nonzero with the given conditions holding. Hence
either two or zero terms are nonzero.  Furthermore, two nonzero terms
have values $\pm e^{du + (d-c) u'}$, with the positive
sign for the first term and the negative for the second or third term,
so they cancel.

In cases where $|\beta|$ or $|\gamma|$ is the largest among $|\alpha|$,
$|\beta|$, and $|\gamma|$, the same argument applies, and the right hand
side is zero.

Next, let us check that $AYBE(y)_{ikm}^{jlm} = 0$ for any $i \neq j, k
\neq l,$ with $i+k = j+l$. We use \eqref{aby}, setting $p = m$.  Set
$\alpha = e_i - e_j$ and $\beta = e_k - e_l$.  It is evident that the first
two terms are each nonzero iff either $-\alpha \peq \beta$ with
$\beta > 0$, or $-\beta \prec \alpha$ with $\alpha > 0$. On the other
hand, the last term is nonzero iff one of these two conditions is
true, with the additional condition that, setting the underscores equal
to $t$, 
either $-\alpha \peq e_m - e_t \prec \beta$, or $-\beta \peq e_t -
e_m \prec \alpha$.  Assuming that all three terms are nonzero,
and using the notational abuse $O(\alpha, \beta) =
O(\text{sign}(\alpha)\alpha, \text{sign}(\beta)\beta)$, we have
\begin{multline}
AYBE(y)_{ikm}^{jlm} = \frac{\text{sign}(\beta)}{1 -
e^{-u-u'}} \bigl( e^{O(\alpha,\beta)u'/n
- O(j,m)(u+u')/n} \\ - e^{-O(\alpha,\beta) u/n - O(k,m)(u+u')/n} 
\bigr) - e^{-O(-\alpha, e_m - e_t) u/n + O(e_m - e_t, \beta) u'/n}.
\end{multline}
Further assuming that $\alpha < 0$ and $-\alpha \prec \beta$, we write
the first two terms of the RHS of \eqref{aby} as
\begin{multline}
\frac{e^{O(j,k) u'/n - O(j,m) (u+u')/n} - e^{-O(j,k) u/n - (O(j,m)/n +
1 - O(j,k)/n)(u+u')}}{1 - e^{-u-u'}} \\ = e^{O(m,k) u'/n - O(j,m) u/n} =
e^{-O(-\alpha, e_m - e_t) u/n + O(e_m - e_t, \beta) u'/n},
\end{multline}
so $AYBE(y)_{ikm}^{jlm} = 0$.  On the other hand, if the third term of
the RHS of \eqref{aby} is zero, and still assuming $\alpha < 0$, then
we can write the first two terms (if nonzero) as
\begin{equation}
\frac{e^{O(j,k)u'/n - [O(j,k) + O(k,m)](u+u')/n} - e^{-O(j,k) u/n - O(k,m)(u+u')/n}}{1 - e^{-u-u'}} = 0,
\end{equation}
so again $AYBE(y)_{ikm}^{jlm} = 0$.  Almost the same thing
happens when $\alpha > 0, -\beta \prec \alpha$.  So, in any case, we
find that $AYBE(y)_{ikm}^{jlm} = 0$.

By the same reasoning, we can see that $AYBE(y)_{ikm}^{jlp} = 0$
whenever either 1) $i = j, k \neq l$, and $m \neq p$ or 2) $k = l, i \neq j$, and $m \neq p$.

Finally, we check that $AYBE(y)_{ikm}^{ikm} = 0$ for all $i, k,$ and
$m$.  We compute:
\begin{multline}
AYBE(y)_{ikm}^{ikm} = \frac{e^{O(i,k)u'/n - O(i,m)(u+u')/n}}{(1 -
e^{u'})(1 - e^{-u-u'})} \\ - \frac{e^{-O(k,m)(u+u')/n - O(i,k)u/n}}{(1
- e^{-u-u'})(1 - e^{-u})} + \frac{e^{-O(i,m)u/n -
O(k,m)u'/n}}{(1-e^{-u})(1-e^{-u'})} \\ = \frac{-e^{-u'+O(i,k)u'/n -
O(i,m)(u+u')/n}(1 - e^{-u}) -e^{-O(k,m)(u+u')/n - O(i,k)u/n}(1 -
e^{-u'})}{(1 - e^{-u})(1 - e^{-u'})(1 - e^{-u-u'})} \\+ \frac{e^{-O(i,m)u/n -
O(k,m)u'/n}(1 - e^{-u-u'})}{(1 - e^{-u})(1 - e^{-u'})(1 - e^{-u-u'})}.
\end{multline}
Let $\delta = 1$ if $i \peq k \peq m$ in the $\tilde
T$-ordering---that is, if $k$ lies between $i$ and $m$ under iteration
of the cyclic permutation $\tilde T$ (or $k = i$ or $m$).  Otherwise,
set $\delta=0$.  Let $\bar \delta$ denote the opposite of $\delta$,
i.e.~$\bar \delta = 1 - \delta$.  Now, we simplify this to:
\begin{multline}
AYBE(y)_{ikm}^{ikm}[(1-e^{-u})(1-e^{-u'})(1-e^{-u-u'})] \\ =
-e^{-O(i,m)u/n - O(k,m) u'/n - u'\delta}(1
- e^{-u}) \\ -e^{-O(i,m)u/n - O(k,m)u'/n - u \bar \delta}(1 -
e^{-u'}) + e^{-O(i,m) u/n -
O(k,m) u'/n}(1 - e^{-u-u'}) \\ =
e^{-O(i,m)u/n - O(k,m) u'/n} [-e^{-u' \delta} + e^{-u' \delta - u} - e^{-u \bar \delta} +
e^{-u \bar \delta - u'} + 1 - e^{-u-u'}] = 0,
\end{multline}
so $AYBE(y)^{ikm}_{ikm} = 0$, independently of $\delta$.  Hence, $y$ satisfies the AYBE.
\end{proof}

\begin{rem} In the preceding proof, the cancellation of terms in the first two
parts of the proof (the ones involving some non-diagonal matrices) is
actually a very special case of the pairing of so-called $T$-quadruples in
\cite{S}.  In \cite{S} these tools are developed much more extensively
to expand the twist from \cite{ESS}, which is an arduous computation.
\end{rem}

\begin{lemma} \label{lau}\cite{P} Let $r$ be a solution of the AYBE
with a Laurent expansion of the form \eqref{laur}.
Then $r$ is uniquely determined by $r_0$.
\end{lemma}
\begin{proof}  
We repeat the computations of \cite{P}.  First, note that, since the
polynomials $u^k, (u')^k,$ and $(u+u')^k$ are linearly independent,
$r_k$ is uniquely determined by $r_0$ and $r_1$ for all $k > 2$.  Now,
from the AYBE for $r$ we obtain the equation
\begin{multline} \label{r01eq}
r_0^{12}(v) r_0^{13}(v+v') - r_0^{23}(v') r_0^{12}(v) + r_0^{13}(v+v')
r_0^{23}(v') \\ = r_1^{12}(v) + r_1^{23}(v') + r_1^{13}(v+v').
\end{multline}
All we have to show is that this equation uniquely determines $r_1$.
Suppose that $r'(u,v)$ is another
AYBE solution with $r'(u,v) = \frac{1 \o 1}{u} + r_0(v) + u r_1'(v) + O(u^2)$.  Then $t = r_1' - r_1$ satisfies
\begin{equation}
t^{12}(v) + t^{13}(v+v') + t^{23}(v') = 0.
\end{equation}
Now, applying $\pr \o id \o id$ to this equation, we obtain $(\pr \o
id) t(v) = 0$ and similarly we obtain $(id \o \pr) t(v) = 0$.  Hence,
$t(v)$ is a scalar meromorphic function satisfying $t(v) + t(v') +
t(v+v') = 0$. Now, for any $k \in \Z \setminus \{0, 1\}$, the elements
$v^k$, $(v')^k$, and $(v+v')^k$ are linearly independent, so when we
write $t$ in terms of its Laurent expansion, we see that the identity
can only be satisfied if $t = a + bv$ for some $a, b \in \C$.  Now the
identity holds iff $a = b = 0$.  Hence, $t(v) = 0$ identically so that
$r_1$ is uniquely given by $r_0$.
\end{proof}

Now, we can complete the
\begin{proof}[Proof of Theorem \ref{mt}, part (2a)]
Uniqueness is a consequence of Lemma \ref{lau}.  By Lemma \ref{ruvc},
$r(u,v)$ indeed has the Laurent expansion \eqref{tle}.  Then, Lemma
\ref{gruvl}, which uses part (3) of the Theorem, reduces our task to
verifying that \eqref{gruv} satisfies the AYBE and the unitarity
condition.  By Lemma \eqref{p12l}, we can assume that $\Phi = 0$,
since the proof of Lemma \ref{gruvl} points out that $\Phi^1 + \Phi^2$
commutes with $r(u,v)$.  By Lemma \ref{nypl}, it suffices only to show
that $y(u)$ given by \eqref{yeq} satisfies the AYBE.  This is proved
in Lemma \ref{ypl}. Hence, the element $r(u,v)$ given by \eqref{afgq}
is a unitary AYBE solution lifting $r_0(v)$, proving part (2a) of
Theorem \ref{mt}.
\end{proof}
\begin{proof}[Proof of Theorem \ref{mt}, part (2b)]
This follows directly from part (2a) and \eqref{rwsp}.
\end{proof}

\subsection{Proof of Theorem \ref{mt}, parts (1a) and (1b)}
\begin{ov}
In this section, we present and exploit condition \eqref{cab}, which
follows from \eqref{r01eq} in Lemma \ref{lau}, in order to prove the
necessity of the associative BD conditions and formula \eqref{stp} for
$s_0$, which is all of (1a) that remains to be proved. The equivalence
of \eqref{stp} and \eqref{tra} is an easy computation, proving part
(1b) and hence the Theorem.
\end{ov}

\begin{lemma}
Suppose that $r(u,v)$ is a solution of the AYBE having a Laurent
expansion of the form \eqref{laur}, where $r_0(v)$ is the classical
$r$-matrix with spectral parameter $r_0(v) = \hat r_{T,s}$
for the BD triple $(\Gamma_1, \Gamma_2, T)$ and matrix $s$.
Then
\begin{equation} \label{cab}
(\pr \o \pr \o \pr) [r_{T,s}^{12} r_{T,s}^{13} - r_{T,s}^{23} r_{T,s}^{12} + r_{T,s}^{13} r_{T,s}^{23}] = 0.
\end{equation}
\end{lemma}
\begin{proof} This follows from \eqref{r01eq} in Lemma \ref{lau},
using the next Lemma \eqref{ovl}.
\end{proof}

\begin{lemma} \label{ovl} Let $r_0(v) = \hat {\tilde r}$
where $r$ satisfies $r + r^{21} = P$.  Then
\begin{multline}
r_0^{12}(v) r_0^{13}(v+v') - r_0^{23}(v')
r_0^{12}(v) + r_0^{13}(v+v') r_0^{23}(v')\\
= \tilde r^{12} \tilde r^{13} - \tilde r^{23} \tilde r^{12} + \tilde r^{13} \tilde r^{23}.
\end{multline}
\end{lemma}
\begin{proof}
Note that $P^{12} \tilde r^{13} = \tilde r^{23} P^{12}$, and similar
relations are all derived from $P t P = t^{21}$.  Substituting $\tilde r^{21} = P - \tilde r^{12}$
six times, we get
$\tilde r^{21} \tilde r^{31} - \tilde r^{32} \tilde r^{21} + \tilde r^{31} \tilde r^{32} = \tilde r^{12} \tilde r^{13} - \tilde r^{23} \tilde r^{12} +
\tilde r^{13} \tilde r^{23}$.  Similarly, we can deduce $\tilde r^{21} \tilde r^{13} - \tilde r^{23} \tilde r^{21} -
\tilde r^{13} \tilde r^{23} = -(\tilde r^{12} \tilde r^{13} - \tilde r^{23} \tilde r^{12} + \tilde r^{13} \tilde r^{23})$ and a
handful of similar identities to expand
\begin{multline}\label{ovrae}
[(1-e^v)(1-e^{v'})
(1-e^{v+v'})][r_0^{12}(v) r_0^{13}(v+v') \\ - r_0^{23}(v')
r_0^{12}(v)
+ r_0^{13}(v+v') r_0^{23}(v')] \\
= (\tilde r^{12} + e^v \tilde r^{21})(\tilde r^{13} + e^{v+v'}
\tilde r^{31})(1-e^{v'}) - (\tilde r^{23} + e^{v'} \tilde r^{32})
(\tilde r^{12} + e^{v} \tilde r^{21})(1 - e^{v+v'})
\\ + 
(\tilde r^{13} + e^{v+v'} \tilde r^{31})
(\tilde r^{23} + e^{v'} \tilde r^{32})(1 - e^v)\\
= (e^{2v+2v'} + e^{2v+v'}+e^{v+2v'}+2 e^{v+v'} + e^{v} + e^{v'} + 1)(\tilde r^{12} \tilde r^{13} - \tilde r^{23} \tilde r^{12} + \tilde r^{13} \tilde r^{23}) \\ =[\tilde r^{12} \tilde r^{13} - \tilde r^{23} \tilde r^{12} + \tilde r^{13} \tilde r^{23}][(1-e^v)(1-e^{v'})
(1-e^{v+v'})],
\end{multline}
proving the Lemma.
\end{proof}
Now, we are in a position to prove the necessity of the BD associativity
conditions given in Definition \ref{atd}:
\begin{lemma}\label{nc1}
\label{orl} The first condition of Definition \ref{atd} is necessary for
an AYBE solution limiting to the CYBE solution to exist.
\end{lemma}
\begin{proof}
Suppose that we are given a Belavin-Drinfeld triple which does not
preserve orientation.  Hence, there exists $i$ and $j$ such that
$T(\alpha_i) = \alpha_{j}$ and
$T(\alpha_{i+1}) = \alpha_{j-1}$.  Now, let $\tilde r$ be the constant solution
of the CYBE corresponding to our Belavin-Drinfeld triple.  Then,
we find that
 $AYBE(\tilde r)_{i+2,j-1,j}^{i,j,j+1}$ $= 1 + 0 + 0 = 1$, so \eqref{cab} is not satisfied.
\end{proof}

\begin{lemma} \label{nc2} The second condition of Definition \ref{atd} is necessary for
the triple to give rise to AYBE solutions.
\end{lemma}
\begin{proof}
We consider the coefficients 
$AYBE(\tilde r)_{i+1,j,k}^{i,j+1,k}$, for $T(\alpha_i) = \alpha_j$. 
We find that
\begin{multline}
AYBE(\tilde r)_{i+1,j,k}^{i,j+1,k} = \tilde r_{i+1, j}^{i, j+1} \tilde r_{i,k}^{i,k} - \tilde r^{j,k}_{j, k} \tilde r_{i+1, j}^{i, j+1} + \tilde r_{i+1, k}^{i, k+1} \tilde r_{j, k+1}^{j+1, k} \\
= [(e_i - e_j) \o e_k] (s + \frac{1}{2} \sum_l e_{ll} \o e_{ll}) - \delta_{ik} = [(e_i - e_j) \o e_k] s - 
\frac{1}{2} (e_i + e_j, e_k).
\end{multline}
In order for $AYBE(\tilde r)$ to be zero modulo scalars, it is necessary
that all of these coefficients are equal for all $k$.  That is, we require
\begin{equation} \label{ntrij}
[(e_i - e_j) \o \alpha] s = \frac{1}{2} (e_i + e_j, \alpha)
\end{equation}
for all roots $\alpha \in \Gamma$. Applying the same work for $AYBE(\tilde r)^{j+1, i, k}_{j, i+1, k}$ we deduce also that
\begin{equation} \label{ntrij2}
[(e_{i+1}- e_{j+1}) \o \alpha] s = \frac{1}{2} (e_{i+1} + e_{j+1}, \alpha)
\end{equation}
for all roots $\alpha$.

Now, provided the first condition of Definition \ref{atd} is satisfied
(which we now know is necessary), we can define a permutation $\tilde
T$ of $\{1, \ldots, n\}$ such that $T(\alpha_i) = \alpha_j$ implies
$\tilde T(i) = j$ and $\tilde T(i+1) = j+1$.  This permutation is
compatible just in the case it is cyclic; we can choose it to be
cyclic iff there is no cycle $(a_1, \ldots, a_k), 1 \leq k < n$, such
that, for each $1 \leq i \leq k$, either $T(\alpha_{a_i}) =
\alpha_{a_{i+1}}$, or $T(\alpha_{a_{i}-1}) = \alpha_{a_{i+1}-1}$
(subscripts of $a$ are given modulo $k$).  Now, in the case that such
a cycle exists, \eqref{ntrij} and \eqref{ntrij2} imply
\begin{equation}
0 = \frac{1}{2} (e_{a_1} + \ldots + e_{a_{k}}, \alpha)
\end{equation}
for any $\alpha$.  This implies that $e_{a_1} + \ldots + e_{a_{k}} = 1$,
so the cycle contains all of $\{1, \ldots, n\}$, contradicting our assumption. 
\end{proof}

\begin{lemma} \label{ays}
Suppose $r(u,v)$ satisfies the AYBE and has a Laurent expansion of the
form \eqref{laur} with $r_0(v) = \hat{\tilde r}$,
where $\tilde r$ is a constant CYBE solution
corresponding to the triple $(\Gamma_1, \Gamma_2, T)$ and $s$.  Write $\tilde
r = a + r_s + s$. Then, for some compatible permutation $\tilde T$,
$s_0 = (\pr \o \pr)s$ satisfies \eqref{tra}.
\end{lemma}

\begin{proof}
Take \eqref{cab} and project to $\h \o \h \o \h$.  Let $t = s +
\frac{1}{2} \sum_i e_{ii} \o e_{ii}$ be the projection of $\tilde r$
to $\h \o \h$.  Define $t'_{ij} = t^{ij}_{ij} - t^{1j}_{1j}
- t^{i1}_{i1}$. Now, \eqref{cab} is equivalent to the condition that
\begin{equation} \label{ch1}
[(e_1 - e_i) \o (e_1 - e_j) \o (e_1 - e_k)] (\tilde r^{12} \tilde r^{13} -
\tilde r^{23} \tilde r^{12} + \tilde r^{13} \tilde r^{23}) = 0
\end{equation}
for all $1 < i, j, k \leq n$. (The same is true if we replace $1$ with any
fixed integer $p$ between $1$ and $n$ and allow $i, j,$ and $k$ to take on any value other than $p$.)  Using the fact that $t^{11}_{11} = \frac{1}{2}$,
we can simplify \eqref{ch1} to
\begin{equation}\label{ch2}
t'_{ij} t'_{ik} - t'_{jk} t'_{ij} + t'_{ik} t'_{jk} = \frac{1}{4}, \quad 1 < i,j,k \leq n.
\end{equation}
Specializing to the case $k = i, i \neq j$, we note that 
$t'_{ij} = - t'_{ji}$, and \eqref{ch2} yields
\begin{equation}
(t'_{ij})^2 = \frac{1}{4}.
\end{equation}
Hence,
\begin{equation}
t'_{ij} = \pm \frac{1}{2}, 1 < i,j \leq n, i \neq j.
\end{equation}
Given the fact that $t'_{ij} = \frac{1}{2}$ and $t'_{jk} =
\frac{1}{2}$ for some distinct $i,j,$ and $k$, \eqref{ch2} implies that
$t'_{ik} = \frac{1}{2}$.  Also, for any distinct $i,j \in \{2,\ldots, n\}$, 
we have $\{t'_{ij}, t'_{ji}\} = \{\frac{1}{2}, -\frac{1}{2}\}$.
Thus, we can obtain a unique total ordering of $\{2, \ldots, n\}$, say the
ordered list $(a_2, \ldots, a_n)$, such that $t'_{a_i a_j} =
\frac{1}{2}$ for all $i < j$.  This is equivalent to a cyclic permutation
of $\{1, \ldots, n\}$ given by $\sigma = (1,a_2,a_3 \ldots,
a_n)$. That is, a cyclic permutation of $(1, \ldots, n)$ is associated
with the ordering of $2, \ldots, n$ obtained by ``cutting off'' $1$.

Evidently the values $t'_{ij}$ completely determine $t$ up to scalars.
We rewrite this in a way which yields \eqref{tra}.
Set $\tilde T = \sigma$. Note that $t_{ii} = \frac{1}{2}$ for $i \neq 1$ and
$t_{11} = - \frac{1}{2}$. Using this, we find
\begin{equation} \label{ntr1}
[(e_i - e_{\tilde T(i)}) \o (e_1 - e_j)] t = t'_{\tilde T(i), j} - 
t'_{ij}
= \delta_{i1} - \delta_{ij} = (e_i, e_1 - e_j).
\end{equation}
Furthermore, we evidently have
\begin{equation} \label{ntr2}
[(e_i - e_{\tilde T(i)}) \o (e_1 - e_j)] \sum_i e_{ii} \o e_{ii} = (e_i - e_{\tilde T(i)}, e_1 - e_j).
\end{equation}
Subtracting one-half of \eqref{ntr2} from \eqref{ntr1}, we get
\begin{equation} \label{ntr}
[(e_i - e_{\tilde T(i)}) \o \alpha] s = \frac{1}{2} 
(e_i + e_{\tilde T(i)}, \alpha)
\end{equation}
for any root $\alpha$.  Letting $P' = P^0 - \frac{1}{n}(1 \o 1) = (\pr
\o \pr) P^0$ denote the projection of $P^0 = \sum_i e_{ii} \o e_{ii}$
to $\mathfrak{g}' \o \mathfrak{g}'$ as in \eqref{tra}, and using
$s_0 = (\pr \o \pr)s$, we can also write \eqref{ntr} as
\begin{equation} \label{ntr3}
[(e_i - e_{\tilde T(i)}) \o 1] s_0 = \frac{1}{2} [(e_i + e_{\tilde T(i)}) \o 1]
P',
\end{equation}
which is exactly \eqref{tra}.

Now, it remains to see that $\tilde T$ is compatible with $T$, that is,
$T(\alpha_i) = \alpha_j$ implies $\tilde T(i) = j$ and $\tilde T(i+1) = j+1$.
To see this, we apply \eqref{ntrij} and \eqref{ntrij2}.
Suppose $T(\alpha_i) = \alpha_j$.
Using the previous work in this lemma,
we know that there is a unique permutation $\tilde T$ such that $s$ satisfies
\eqref{tra} for all roots $\alpha$.  
In particular, \eqref{tra} implies that
\begin{equation} \label{ijtt}
[(e_i - e_j) \o \alpha] s = \frac{1}{2}(e_i + e_j, \alpha) + \sum_{k : \ 0 < O(i,k) < O(i,j)} (e_k, \alpha), \forall \alpha \in \Gamma.
\end{equation}
Equating this with the right-hand side of \eqref{ntrij}, we conclude
that $j = \tilde T(i)$.  Also, \eqref{ijtt} continues to hold
replacing $i$ and $j$ with $i+1$ and $j+1$, respectively.  Comparing
this with \eqref{ntrij2}, we also find that $\tilde T(i+1) =
j+1$. This completes the proof.
\end{proof}

\begin{lemma} \label{stptra}
The condition \eqref{stp} is equivalent to the condition \eqref{tra}.
\end{lemma}
\begin{proof}
Let $\h_0 = \h \cap \mathfrak{g}'$ be the space of traceless diagonal
matrices.  Since \eqref{tra} uniquely determines $s_0 \in \h_0 \wedge \h_0$,
it suffices to show that the element
$s_0$ given by \eqref{stp} satisfies \eqref{tra}.
This is verified as follows, letting $s_0$ be given by \eqref{stp}:
\begin{multline}
[(e_i - e_{\tilde T(i)}) \o 1] s_0 = \sum_{j \notin \{i, \tilde T(i)\}}
- \frac{1}{n} e_{jj} + \frac{n-2}{2n} (e_{ii} + e_{\tilde T(i), \tilde
T(i)}) \\ = \frac{1}{2} (e_i + e_{\tilde T(i)} \o 1) P',
\end{multline}
as desired.  
\end{proof}

\begin{proof}[Proof of Theorem \ref{mt}, parts (1a) and (1b)]
Part (2a) proves sufficiency of the conditions, since the $r(u,v)$ given
by \ref{afgq} lifts $r_0(v) = r_{T,s}$ and is a unitary AYBE solution (and
it satisfies the BD associativity and $s_0$-conditions).
For the ``only-if,'' or necessity, Lemma \ref{ays} proves that, given
any AYBE solution $r(u,v)$ with a Laurent expansion as in
\eqref{laur}, for $r_0(v) = r_{T,s}$, $s$ must satisfy \eqref{tra} for
some unique compatible permutation $\tilde T$. In particular, this
implies that the BD triple is associative, which alternately follows
from Lemmas \ref{nc1} and \ref{nc2}.  This completes the proof of (1a).

By Lemma \ref{stptra}, part (1b) follows.
\end{proof}

This completes the proof of Theorem \ref{mt}.

\bibliography{trigaybe}
\bibliographystyle{amsalpha}

\noindent
Address: 45 rue d'Ulm; 75005 PARIS; France.

\noindent
Email: {\sl schedler@post.harvard.edu}

\end{document}